\documentclass[12pt,reqno]{amsart}
\usepackage{hyperref}
\usepackage{enumerate,fullpage}
\usepackage{amssymb,amsmath,graphicx,amsthm}
\usepackage{color}
\usepackage{soul}

\def\B{\4B}

\def\Om{\Omega}

\def\no#1{\|#1\|}

\numberwithin{equation}{section}
\def\db{\bar\partial}
\def\db*{\bar\partial^*}
\def\T{\text}
\def\simgeq{\gtrsim}
\def\simleq{\lesssim}
\def\1#1{\overline{#1}}
\def\2#1{\widetilde{#1}}
\def\3#1{\widehat{#1}}
\def\4#1{\mathbb{#1}}

\def\5#1{\frak{#1}}
\def\6#1{{\mathcal{#1}}}
\def\C{{\4C}}

\def\B{\4B}

\def\di{\partial}
\def\dib{\bar\partial}

\def\Re{{\sf Re}\,}

\def\phi{\varphi}

\def\Om{\Omega}

\newtheorem{Thm}{Theorem}[section]
\newtheorem{Cor}[Thm]{Corollary}
\newtheorem{Pro}[Thm]{Proposition}
\newtheorem{Lem}[Thm]{Lemma}
\theoremstyle{definition}\newtheorem{Def}[Thm]{Definition}
\theoremstyle{remark}
\newtheorem{Rem}[Thm]{Remark}
\newtheorem{Exa}[Thm]{Example}

\def\Label#1{\label{#1}}
\def\bl{\begin{Lem}}
\def\el{\end{Lem}}
\def\bp{\begin{Pro}}
\def\ep{\end{Pro}}
\def\bt{\begin{Thm}}
\def\et{\end{Thm}}
\def\bc{\begin{Cor}}
\def\ec{\end{Cor}}
\def\bd{\begin{Def}}
\def\ed{\end{Def}}
\def\br{\begin{Rem}}
\def\er{\end{Rem}}
\def\be{\begin{Exa}}
\def\ee{\end{Exa}}
\def\bpf{\begin{proof}}
\def\epf{\end{proof}}
\def\ben{\begin{enumerate}}
\def\een{\end{enumerate}}

\def\1alpha{[\frac1\alpha]}
\def\T{\text}

\def\C{{\Bbb C}}

\def\Real{\Re}

\begin{document}
\title[H\"older regularity...]{H\"older regularity of the solution to the complex Monge-Amp\`ere equation with $L^p$ density}
\author[L.~Baracco, T.V.~Khanh, S.~Pinton and G.~Zampieri]{Luca Baracco, Tran Vu Khanh, Stefano Pinton and  Giuseppe Zampieri}
\address {Luca Baracco, Stefano Pinton, Giuseppe Zampieri}
\address{Dipartimento di Matematica, Universit\`a di Padova, via 
Trieste 63, 35121 Padova, Italy}
\email{baracco@math.unipd.it,pinton@math.unipd.it}

	\address{Tran Vu Khanh}
	\address{School of Mathematics and Applied Statistics, University of Wollongong, NSW, Australia,  2522}
	\email{tkhanh@uow.edu.au}
\thanks{The research of T.~V.~Khanh was supported by the Australian Research Council
	DE160100173.}
\maketitle

\abstract
On a smooth domain $\Om\subset\subset\C^n$, 
we consider the Dirichlet problem for the complex Monge-Amp\`ere equation $((dd^cu)^n=fdV,\,u|_{b\Om}\equiv\phi)$. We state the H\"older regularity of the solution $u$ when the boundary value $\phi$ is H\"older continuous and the density $f$ is only $L^p$, $p>1$. Note that in former literature (Guedj-Kolodziej-Zeriahi) the weakness of the assumption $f\in L^p$ was balanced  by taking $\phi\in C^{1,1}$ (in addition to assuming $\Om$ strongly pseudoconvex). 

\noindent
MSC: 32U05, 32U40, 53C55 
\endabstract
%\endtopmatter
\section{Introduction}
For a bounded pseudoconvex domain $\Om\subset\subset \C^n$, the Dirichlet problem for the Monge-Amp\`ere equation consists in
\begin{equation}
\Label{1.1}
\begin{cases}
(dd^cu)^n=fdV&\T{in }\Om,
\\
u=\phi&\T{on }b\Om.
\end{cases}
\end{equation}
In our discussion we take a density $0\le f\in L^p(\Om)$, $ \,1< p\le+\infty$, a boundary datum $\phi$ in some H\"older class and look for a plurisubharmonic solution $u\in  C^\beta(\bar\Om)$ for a certain $\beta$. 
Sometimes we use the notation $MA(\phi,f)$ for the problem \eqref{1.1} and $u(\Om,\phi,f)$ for its solution. This problem has been extensively investigated in recent years under the assumption that $\Om$ is strongly pseudoconvex.  Bremermann \cite{B59}, Walsh \cite{W68} and Bedford-Taylor\cite{BT76} show that there is a solution $u\in C^0(\bar\Om)$ if $\phi\in C^0(b\Om)$ and $f\in C^0(\bar\Om)$. By the well known ``comparison principle"' (cf. Kolodziej \cite{K05}), the solution is unique; what matters   is to prove the H\"older continuity of this $C^0$-solution. In this direction, in \cite{BT76} is proved that $u\in C^{\frac{\alpha}{2}}(\bar\Om)$ if $\phi\in C^{\alpha}(b\Om)$, $f^{\frac1n}\in C^{\frac\alpha2}(\bar\Om)$. A recent interest has been dedicated to the case when $\Om$ is no longer strongly pseudoconvex but has a certain ``finite type" $m$. Li proves in \cite{L04} that $u\in C^{\frac\alpha m}(\bar\Om)$ if $\phi\in C^{\alpha}(b\Om)$ and $f^{\frac1n}\in C^{\frac\alpha m}(\bar{\Om})$. Ha and Khanh in \cite{HK14} get the same conclusion with a more geometric notion of finite type (cf. \eqref{1.2} below) and have also a generalization for the infinite type. Coming back to the case of $\Om$ strongly pseudoconvex, Caffarelli, Kohn and Nirenberg prove in \cite{CKNS85} that $u\in C^\infty(\bar{\Om})$, for $\phi\in C^\infty(b\Om)$ and $f\in C^\infty(\bar{\Om})$, in case $f>0$ in $\bar\Om$. Lowering the smoothness of $f$ gives the problem additional difficulty. Guedj, Kolodziej and Zeriahi prove in \cite{GKZ08} that if $f\in L^p(\Om)$ with $p>1$ and $\phi\in C^{1,1}(b\Om)$ then $u\in C^{\gamma}(\bar\Om)$ for any $\gamma<\gamma_p:=\frac2{qn+1}$ where $\frac1q+\frac1p=1$. Recently, Charabati has obtained in  \cite{Ch15} that $u\in C^{\frac\gamma{2}}(\bar\Om)$ for the same datum as in \cite{GKZ08} on a bounded strongly hyperconvex Lipschitz domain i.e. on a domain for which there exists a Lipschitz  plurisubharmonic defining function $\rho$ such that $(dd^c\rho)^n\ge cdV$. Our purpose is twofold: to lower the regularity of $\phi$ and to allow a (geometric) finite type $m$ for $\Om$ with some $m\ge 2$. What we get is that  if $f\in L^p(\Om)$ with $p>1$ and $\phi\in C^{\alpha}(b\Om)$ with $0<\alpha\le 2$ then $u\in C^{\frac{\alpha}{m}}(\bar{\Om})$ if $\alpha<\gamma_p$ otherwise $u\in C^\frac{\gamma}{m}(\bar{\Om})$ for any $\gamma<\gamma_p$. To go into the detail of our geometric setting we consider a submanifold $S\subset b\Om$  of CR dimension $0$. Let $d_S$ be the distance to $S$ and $(L_{b\Om})$ be the Levi form of $b\Om$. We assume that $b\Om$ has finite type $m$ along $S$  in the sense that  
\begin{equation}
\Label{1.2}
L_{b\Om}\simgeq d_S^{m-2}.
\end{equation}
To convert \eqref{1.2} into a suitable property for our use, we need two basic results. First, from Khanh and Zampieri \cite{KZ10}, we know that \eqref{1.2} implies the potential-theoretic  ``$t^{\frac1m}$-property". By \cite{K13} and \cite{HK14} this implies in turn that there is an exhaustion function $\rho$ which defines $\Om $ by $\rho<0$ such that
\begin{equation}
\Label{1.3}
i\di\dib\rho\ge \T{Id}\,\,\T{ in $\Om$},
\quad \rho\in C^{\frac2m}(\bar{\Om}).
\end{equation}
\br
According to Catlin \cite{C87}, if $\Om$ has finite D'Angelo type $D$, then it has the ``$t^{\frac1m}$-property" for $\frac1m:=D^{-n^2D^{n^2}}$; again, this implies the existence of the exhaustion $\rho\in C^{2D^{-n^2D^{n^2}}}(\bar{\Om})$  with $i\di\dib\rho\ge \T{Id}$ in $\Om$.
\er
It is \eqref{1.3} the property which rules many passages of this paper.
Here is our result
\bt
\Label{t1.1}
Let $\Om\subset\subset\C^n$ be a $C^2$-smooth pseudoconvex domain of finite type $m$ with $m\ge 2$ in the sense of \eqref{1.2} and let $\phi\in C^{\alpha}(b\Om)$ with $0<\alpha\le2$ and $f\in L^p(\Om)$ with $p>1$. Then the unique solution $u$ to $MA(\Om,\phi,f)$  is in $C^{\min(\frac{\alpha}{m},\frac\gamma{m})}(\bar{\Om})$ with $\gamma<\gamma_p$ where $\gamma_p:=\frac2{qn+1}$ and $\frac 1 p + \frac 1 q$=1.
\et
The proof follows in Section~\ref{s3}. %\red{Note that, in particular, if $\Om$
%is strongly pseudoconvex and $\alpha=2$, then $u\in C^\gamma(\bar{\Om})$. It coincides with the result of \cite{GKZ08}.}\\

Throughout the paper we use $\lesssim$ and $\gtrsim$ to denote an estimate up to a positive constant and $\sim$ for the combination of $\lesssim$ and $\gtrsim$. Finally, the indices $m$, $\alpha$,  $p$, $\gamma$ and $\gamma_p$ only take ranges as in Theorem~\ref{t1.1}.  

\vskip0.3cm
\noindent
{\it Acknowledgments.}\hskip0.3cm The paper is largely inspired by  S. Kolodziej on the stream of research initiated by Bedford and Taylor in 1976. The authors are deeply indebted to S. Kolodziej for his important advice during private communications. The authors are also grateful to the referee for his helpful comments.
\section{H\"older regularity of a subsolution}
\Label{s2}
We say that $v\in C^0(\bar{\Om})$ is a subsolution to $MA(\Om,\phi,f)$ if $v$ is plurisubharmonic, $v|_{b\Om}=\phi$ and $(dd^cv)^n\ge f$ in $\Om$.
\bp
\Label{p2.2}
Let $\rho$ satisfy \eqref{1.3}. Then there is a subsolution $v\in C^0(\bar\Om)$ to $MA(\
\Om,\phi,f)$ for $\phi\in C^0(b{\Om})$ and $f\in L^p(\Om)$.
\ep
\bpf
For a large ball $\B$ containing $\Om$, we define
\begin{equation*}
\tilde f(z):=\begin{cases} f(z)\quad\T{if $z\in\Om$},
\\
0\quad\T{ if $z\in\B\setminus\Om$}.
\end{cases}
\end{equation*}
We consider the solutions
\begin{equation*}
\begin{cases}
u_1=u(\B,0,\tilde f)\in C^0(\bar{\B})\quad \T{by Kolodziej on the ball $\B$ (strongly pseudoconvex) \cite{K98}},
\\
u_2=u(\Om,(-u_1)|_{b\Om},0)\in C^0(\bar{\Om})\quad\T{by Blocki \cite{B96}}.
\end{cases}
\end{equation*}
Taking summation $u_1+u_2$ we have a subsolution to $MA(\Om,0,f)$ in $C^0(\bar{\Om})$. Using the solution $u(\Om,\phi,0)\in C^0(\bar \Om)$ provided by \cite{B96} and putting
$$
v=u_1+u_2+u(\Om,\phi,0),
$$
we get the desired subsolution.

\epf

We change a little our setting and take 
$\phi\in C^{\alpha}(b\Om)$ and $f\in L^\infty(\Om)$. If $\zeta$ is a general point of $b\Om$ we set
\begin{equation}
\Label{2.1}
v_\zeta(z):=\begin{cases}
\phi(\zeta)-c[-\rho(z)+|z-\zeta|^2]^{\frac{\alpha}{2} }\quad\T{ if $0<\alpha\le1$},
\\
\phi(\zeta)-\sum_j2\Real \frac{\di\phi}{\di z_j}(\zeta)(z_j-\zeta_j)-c[-\rho(z)+|z-\zeta|^2]^{\frac\alpha2}\quad\T{ if $1<\alpha\le2$}.
\end{cases}
\end{equation}
If there is an exhaustion function $\rho\in C^{\frac2m}(\bar{\Om})$ such that $i\di\dib\rho\ge\T{Id}$ in $\Om$ then we can find $c$, independent of $\zeta$ and only depending on $\no{\phi}_{C^\alpha(\bar{\Om})}$ and $\no{f}_{L^\infty(\Om)}$ such that (cf. \cite{HK14,L04})
\begin{equation}
\Label{2.2}
\begin{cases}
v_\zeta(z)\le\phi(z)\quad\T{ if $z\in b\Om$},
\\
v_\zeta(\zeta)=\phi(\zeta),
\\
(dd^cv_\zeta)^n\ge f \quad\T{in } \Om,
\\
v_\zeta\in C^{\frac\alpha{m}}(\bar{\Om}).
\end{cases}
\end{equation}
Using the family $\{v_\zeta\}_{\zeta\in b\Om}$ it is readily seen (cf. \cite{HK14,L04}) that for any plurisubharmonic $C^0(\bar{\Om})$ solution  to MA we have $u(\Om,\phi,f)\in C^{\frac\alpha{m}}(\bar{\Om})$ for $\phi\in C^{\alpha}(b\Om)$ and $f^{\frac1n}\in C^{\frac\alpha{m}}(\bar{\Om})$; in particular, $u(\Om,\phi,0)\in C^{\frac\alpha{m}}(\bar{\Om})$ for $\phi\in C^{\alpha}(b\Om)$. We lower the smoothness of $f$. We start from
\bp
\Label{p2.1}
Let $\rho$ satisfy \eqref{1.3}.  Then there is a subsolution $v\in C^\frac{\alpha}{m}(\bar\Om)$ to $MA(\
\Om,\phi,f)$ for $\phi\in C^\alpha(b{\Om})$ and $f\in L^\infty(\Om)$.
\ep

%We start from a subsolution $v\in C^{\frac1m}$ when $\phi=0$. We set $v=c\rho$; thus $v|_{b\Om}\equiv0$ and $v\in C^{\frac1m}$. We select $K\subset\subset \Om$, take $c\ge \underset{\Om}\sup\,|f|,\,\,c\ge\underset K \sup\,\frac{|
\bpf
We consider the solution $u(\Om,\phi,0)\in C^{\frac\alpha m}(\bar{\Om})$ by \cite{L04} and \cite{HK14} and define
$$
v=u(\Om,\phi,0)+c\rho.
$$ 
For $c\simgeq \no{f}_{L^\infty(\Om)}^\frac{1}{n}$, $v$ is a subsolution.

\epf

We now take $f\in L^p(\Om)$.
\bp
\Label{p2.3}
Let $\rho$ satisfy \eqref{1.3}.Then there is a subsolution $v\in C^{\min(\frac{\alpha}{m},\frac{\gamma}{m})}(\bar\Om)$ to $MA(\
\Om,\phi,f)$ for $\phi\in C^\alpha(b{\Om})$ and $f\in L^p(\Om)$.
\ep

\bpf
We define $\B$ and $\tilde f$ as in the proof of Proposition~\ref{p2.2}. Since $\tilde f$ is bounded near the boundary, we consider the solutions
\begin{equation*}
\begin{cases}
u_1=u(\B,0,\tilde f)\in C^{\gamma}(\bar{\B})\quad \T{by \cite{GKZ08}},
\\
u_2=u(\Om,(-u_1)|_{b\Om},0)\in C^{\frac {\gamma} m}(\bar{\Om})\quad\T{by \cite{L04} and \cite{HK14}}.
\end{cases}
\end{equation*}
Taking the solution $u(\Om,\phi,0)\in C^{\frac\alpha{m}}(\bar{\Om})$ (cf. \cite{HK14,L04}) and taking summation $v=u_1+u_2+u(\Om,\phi,0)$ we have the conclusion.

\epf

\section{H\"older regularity of the solution - Proof of Thorem~\ref{t1.1}}
\Label{s3}
We recall a crucial fact from \cite{K98}. For a general domain, not necessarily strongly pseudoconvex, the existence of $u(\Om,\phi,0)\in C^0(\bar{\Om})$ (which turns out to be equivalent to the existence of a maximal function with boundary datum $\phi$), in addition to the existence of a subsolution $v\in C^0(\bar{\Om})$ for $\phi\in C^0(b\Om)$ and $f\in L^p(\Om)$,  implies the existence of a solution $u(\Om,\phi,f)\in L^\infty(\Om)$. %Furthermore, by \cite{K102}, this solution is in fact in $C^0$. 
In particular, 
\bt (Kolodziej \cite{K98})
\Label{t3.0}
Assume $\Om$ is defined by $\rho<0$ for $\rho\in C^0(\bar\Om)$ such that $i\di\dib\rho\ge\T{Id}$ in $\Om$. Then for any $\phi\in C^0(b\Om)$, $f\in L^p(\Om)$ there is a (unique) plurisubharmonic solution $u(\Om,\phi,f)\in L^\infty(\Om)$.
\et
\bpf
By the property of $\rho$, which implies b-regularity, there is a solution for continuous data, in particular for $f=0$, that is $u(\Om,\phi,0)$ (cf. \cite{B96}); thus there is a maximal function for the given boundary data. Again by the property of $\rho$, there is a subsolution for $\phi\in C^0(b\Om)$, $f\in L^p(\Om)$ (Proposition~\ref{p2.2} above).  
Then by  \cite{K98} Thm. C p. 97 (3 lines after the statement) there is a solution in $L^\infty(\Om)$.
%and by \cite{K02} this solution is in fact in $C^0$.

\epf
\br
\Label{r3.1}
The solution $u(\Om,\phi,f)$ for $\phi\in C^0(b\Om),\,f\in L^p(\Om)$ is in fact in $C^0(\bar\Om)$ by Kolodziej \cite{K02}. Note that the paper makes the general assumption of pseudoconvexity of $\Om$ but this is needless for this specific conclusion. This is confirmed by private communication with the author.
\er
We assume from now $i\di\dib\rho\ge\T{Id}$ in $\Om$ for $\rho\in C^{\frac2m}(\bar{\Om})$.
According to Proposition~\ref{p2.3} above, when we take a smoother boundary datum $\phi\in C^{\alpha}(b\Om)$, there is a subsolution $v\in C^{\min(\frac{\alpha}{m},\frac{\gamma}{m})}(\bar\Om)$ for $f\in L^p$. What follows  is dedicated to show that, in this situation, the $L^\infty$ plurisubharmonic solution $u(\Om,\phi,f)$ is in fact in $C^{\min(\frac{\alpha}{m},\frac{\gamma}{m})}(\bar\Om)$.\\

Let $w:=u(\Omega,\phi,0)\in C^{\frac{\alpha}{m}}(\bar\Om)$ (cf. \cite{HK14,L04}); comparison principle yields at once
\begin{equation}
\Label{3.1}
v\le u(\phi,f)\le w.
\end{equation}
By \eqref{3.1} and by the $C^{\min(\frac{\alpha}{m},\frac{\gamma}{m})}$ regularity of $v$ and $w$ we get
$$
|u(z)-u(\zeta)|\simleq |z-\zeta|^{\min(\frac{\alpha}{m},\frac{\gamma}{m})},\quad z\in \bar\Om,\,\,\zeta\in b\Om,
$$
and therefore for $\delta$ suitably small
\begin{equation}
\Label{3.0}
|u(z)-u(z')|\simleq\delta^{\min(\frac{\alpha}{m},\frac{\gamma}{m})},\quad z,\,z'\in\Om\setminus \Om_\delta\T{ and }|z-z'|<\delta
\end{equation}
 where $\Omega_{\delta}:=\{ z\in\C^n:\, r(z)<-\delta \}$ and $r$ is a $C^2$ defining function for $\Omega$ with $|\nabla r|=1$  in a neighborhood of $b\Omega$. We have to prove that \eqref{3.0} also holds for $z,\,z'\in\Om_\delta$. We use the notation
\begin{equation}
\Label{**}
\begin{cases}
u_{\frac\delta2}:=\underset{|\zeta|<{\frac\delta2}}\sup u(z+\zeta),\quad z\in\bar\Om_\delta,
\\
\tilde u_{\frac\delta2}:=\frac1{\sigma_{2n-1}{\big(\frac\delta2}\big)^{2n-1}}\int_{b\B(z,{\frac\delta2})}u(\zeta)dS(\zeta),\quad z\in\bar\Om_\delta,
\end{cases}
\end{equation}
where ${\sigma_{2n-1}{\big(\frac\delta2}\big)^{2n-1}}=Vol(b\B(z,\frac\delta2))$. 
It is a classical consequence of Riesz Theorem that for a general plurisubharmonic function $u$, not necessarily $C^2$, there is well defined $\Delta u$ in the space of positive Borel measures. We use the notation $\no{\Delta u}^\Om$ for the total mass of $\Delta u $ on $\Om$. 
\bt
\Label{t3.1}
Let $0<\epsilon<1$. We have
\begin{equation}
\Label{3.2}
\no{\tilde u_{\frac\delta2}-u}_{L^1(\Om_\delta)}\simleq \delta^{1-\epsilon}\no{(-r)^{1+\epsilon} \Delta u}^{\Om_{\frac\delta2}}.
\end{equation}
\et

\noindent
%{\it Proof.}\hskip0.2cm
\bpf
The proof is inspired by \cite{GKZ08} Lemma~4.3; the novelty here consists in replacing $\delta^2$ by $\delta^{1-\epsilon}(-r)^{1+\epsilon}.$ 
We start from
\begin{equation}
\Label{3.10}
\begin{split}
\tilde u_{\frac\delta2}(z)-u(z)&\sim \frac1{\delta^{2n-1}}\int_{b\B(0,\frac\delta2)}(u(z+\xi)-u(z))dS(\xi)
\\
&\sim \frac1{\delta^{2n-2}}\int_{b\B(0,\frac\delta2)}dS(\xi)\int_0^1\nabla u(z+s\xi)\cdot \frac\xi\delta ds
\\
&\underset{\T{divergence thm.}}=\frac1{\delta^{2n-2}}\int_0^1sds\int_{\B(0,\frac\delta2)}\Delta u(z+s\xi)
\\
&\underset{s\xi=\zeta,\,s\delta=t}\sim\frac1{\delta^{2n-2}}\int_0^{\frac\delta2}\frac{t}{\delta^2}\frac{t^{-2n}}{\delta^{-2n}}dt\int_{\B(0,t)}\Delta u(z+\zeta).
\end{split}
\end{equation}

%%%%%%
We denote by $\tau_\zeta$ the translation by $\zeta$ and observe that $\tau_\zeta\Om_\delta\subset\Om_{\frac\delta2}\subset\subset\Om$ for any $\zeta\in \B(0,t)$. Observing that the positive measure $\Delta u$ has finite mass on compact subsets of  $\Om$, in particular on $\Om_{\frac\delta2}$, we get, for $t<\frac\delta2$
\begin{equation}
\Label{supernova}
\int_{\Om_\delta}dV(z)\int_{\B(0,t)}\Delta u(z+\zeta)\simleq t^{2n}\int_{\Om_{\frac\delta2}}\Delta u(z).
\end{equation}
We now perform integration $\int_{\Om_\delta}\cdot dV(z)$ in both sides of \eqref{3.10}, apply \eqref{supernova} and end up with
\begin{equation}
\begin{split}
\Label{3.11}
\int_{\Om_\delta}(\tilde u_{\frac\delta 2}-u)(z)\,dV(z)&
\simleq \int_0^{\frac\delta2}t^{-2n+1}t^{2n}dt\int_{\Om_{\frac\delta2}}\Delta u\\
&\simleq \int_0^{\frac\delta2}t\delta^{-(1+\epsilon)}dt\int_{\Om_{\frac\delta2}}(-r)^{1+\epsilon}\Delta u
\\
&\sim \delta^{1-\epsilon}\no{(-r)^{1+\epsilon} \Delta u}^{\Om_{\frac\delta2}}.
\end{split}
\end{equation}

\epf

At this point, the problem is to prove the boundedness of $\no{(-r)^{1+\epsilon}\Delta u}^{\Om_{\frac\delta2}}$ uniformly in $\delta$. This holds (cf. Theorem~\ref{t3.2} below) because of the presence of the factor $(-r)^{1+\epsilon}$. In absence of this factor, one should suppose from the beginning that $\Delta u$ has finite total mass on $\Om$; in turn, this would be a consequence of the hypothesis $\phi\in C^{1,1}$ (cf. \cite{GKZ08}). 
\bt
\Label{t3.2}
We have
\begin{equation}
\no{(-r)^{1+\epsilon}\Delta u}^{\Om}\simleq \no{(-r)^{-1+\epsilon}u}_{L^1(\Om)}.
\end{equation}
\et
\bpf
%It is not restrictive to assume that $u$ has  support in $\Om'\subset\subset\Om$; for this, we just have to multiply $u$ by a cut-off which is $1$ on $\Om_{\frac\delta2}$.
We take a system of smooth cut-off functions $\chi_\nu(|z|)\in C^\infty_c(\B^{2n}(0,\frac1\nu))$, $\no{\chi_\nu}_{L^1}\equiv1$, $\frac1\nu\to0$, and regularize
$$
u_\nu:=\int_\Om u(\tau)\chi_\nu(|z-\tau|)dV(\tau).
$$
 The $u_\nu$'s belong to $C^\infty(\Om)$, converge to $u$ on $\Om$, and satisfy
\begin{equation}
\Label{3.16}
\begin{cases}
\underset{\Om_{\frac1\nu}}\sup\,|\nabla u_\nu|=\underset{\Om_{\frac1\nu}}\sup\,|\nabla(u*\chi_\nu)|\le\nu\no{u}_{L^1(\Om)}
\\
\underset{\Om_{\frac1\nu}}\sup\,u_\nu\le c\quad\T{independent of $\nu$.}
\end{cases}
\end{equation}
Now that the $u_\nu$'s are regular, the $\Delta u_\nu$'s are well defined functions and hence we use the notation $\Delta u_\nu dV$ for the associated measures.
We have 
\begin{equation}
\Label{3.17}
\begin{split}
\int_{\Om_{\frac1\nu}}&(-r)^{1+\epsilon}\Delta u_\nu dV(z)=\int_{\Om_{\frac1\nu}}\T{div}((-r)^{1+\epsilon}\nabla u_\nu)dV(z)+(1+\epsilon)\int_{\Om_{\frac1\nu}}(-r)^\epsilon\nabla r\cdot\nabla u_\nu dV(z)
\\
&\underset{\T{Stokes}}=\int_{b\Om_{\frac1\nu}}(-r)^{1+\epsilon}\nabla r\cdot\nabla u_\nu dS^{2n-1}(z)+(1+\epsilon)\int_{\Om_{\frac1\nu}}(-r)^\epsilon\nabla r\cdot\nabla u_\nu dV(z)
\\
&=\int_{b\Om_{\frac1\nu}}(-r)^{1+\epsilon}\nabla r\cdot\nabla u_\nu dS^{2n-1}(z)+(1+\epsilon)\int_{\Om_{\frac1\nu}}\T{div}((-r)^\epsilon(\nabla r\,u_\nu))dV(z)
\\&\hskip0.3cm+\epsilon(1+\epsilon)\int_{\Om_{\frac1\nu}}(-r)^{\epsilon-1}\nabla r\cdot \nabla r u_\nu dV(z)-(1+\epsilon)\int_{\Om_{\frac1\nu}}(-r)^\epsilon\Delta ru_\nu dV(z)
\\
&\underset{\T{Stokes}}=
\int_{b\Om_{\frac1\nu}}(-r)^{1+\epsilon}\nabla r\cdot\nabla u_\nu dS^{2n-1}(z)+(1+\epsilon)\int_{b\Om_{\frac1\nu}}(-r)^\epsilon\nabla r\cdot\nabla r\,u_\nu dV(z)
\\
&\hskip0.3cm+\epsilon(1+\epsilon)\int_{\Om_{\frac1\nu}}(-r)^{\epsilon-1}\nabla r\cdot \nabla r u_\nu dV(z)-(1+\epsilon)\int_{\Om_{\frac1\nu}}(-r)^\epsilon\Delta ru_\nu dV(z)
\\
&\underset{\T{\eqref{3.16}}}\simleq O(\nu^{-\epsilon})+(1+\epsilon)O(\nu^{-\epsilon})+\int_{\Om_{\frac1\nu}}(-r)^{\epsilon-1}|u_\nu|dV(z)+\int_{\Om_{\frac1\nu}}(-r)^\epsilon |u_\nu|dV(z)
\\
&\simleq O(\nu^{-\epsilon})+ \no{(-r)^{-1+\epsilon}u}_{L^1(\Om)}.
\end{split}
\end{equation}
On the other hand, since $u$ is plurisubharmonic, then $\Delta u$ is a measure on $\Om$ and $\Delta u_\nu dV\underset{\T{weakly}}\to \Delta u.$
The conclusion follows from the following elementary Lemma 
\bl
\Label{l3.1}
Assume $\Delta u_\nu\ge0$ and 
\begin{equation*}
\begin{cases}
\int_{\Om_{\frac1\nu}}(-r)^{1+\epsilon}\Delta u_\nu dV \T{ are bounded}
\\
\Delta u_\nu dV\underset{\T{weakly}}\to \Delta u\T.
\end{cases}
\end{equation*}
Then
$$
\T{$ \int_\Om(-r)^{1+\epsilon}\Delta u$ is  bounded.}
$$
\el
The proof is just a consequence of the dominated convergence theorem for the sequence $(-r)^{1+\epsilon}\psi_\nu\Delta u_\nu dV\to (-r)^{1+\epsilon}\Delta u$ where $\psi_\nu$ are the characteristic functions of the sets $\Om_{\frac1\nu}$. With Lemma~\ref{l3.1} in our hands, we 
get the conclusion of the proof of Theorem~\ref{t3.2}.

\epf

 To end the proof of Theorem \ref{t1.1} we shall need the stability estimate (Theorem $(1.1)$ in \cite{GKZ08})
\bt\label{tse}
Fix $0\leq f\in L^p(\Omega),\, p>1$. Let $U, \, W$ be two bounded plurisubharmonic functions in $\Omega$ such that $(dd^c U)^n=fdV$ in $\Omega$ and let $U\geq W$ on $\partial\Omega$. Fix $s\geq 1$ and $0\leq \eta< \frac{1}{nq+s}$, $\frac 1p+\frac 1q=1$. Then there exists a uniform constant $C=C(\eta, \|f\|_{L^p(\Omega)})>0$ such that 
$$
\sup_\Om(W-U)\leq C\no{(W-U)_+}^\eta_{L^s(\Om)},
$$
where $(W-U)_+:=\max(W-U,0)$.
\et

\bpf[End of Proof of Theorem~\ref{t1.1}.]
Again, we follow the guidelines of \cite{GKZ08}.
Along with $\tilde u_\delta$ defined by \eqref{**} we introduce $\hat u_\delta:=\frac1{\sigma_{2n}{(\delta})^{2n}}\int_{\B(z,{\delta})}u(\zeta)dV(\zeta)$, $\quad z\in\Om_\delta$. We recall that Lemma 4.2 of \cite{GKZ08} states the equivalence between
\begin{equation}
\Label{a}
\underset {\Om_\delta}\sup(u_\frac\delta2-u)\lesssim \delta^{\min(\frac{\alpha}{m},\frac{\gamma}{m})}
\end{equation}
and
\begin{equation}
\Label{b}
\underset {\Om_\delta}\sup(\hat{u}_\frac\delta2-u)\lesssim \delta^{\min(\frac{\alpha}{m},\frac{\gamma}{m})}
\end{equation}
On the other hand, on account of the obvious inequalities
$$
\hat u_\delta\le \tilde u_\delta\le u_\delta,
$$ we see that whatever of \eqref{a} and \eqref{b} is equivalent to
\begin{equation}
\Label{*}
\underset {\Om_\delta}\sup(\tilde{u}_\frac\delta2-u)\lesssim \delta^{\min(\frac{\alpha}{m},\frac{\gamma}{m})}.
\end{equation}
We have thus to prove \eqref{*}. To see it, we remark that 
\begin{equation}
\Label{hypernova}
\begin{split}
\no{\tilde u_{\frac\delta2}-u}_{L^1(\Om_\delta)}&\underset{\T{Theorem~\ref{t3.1}}}\lesssim\delta^{1-\epsilon}\no{(-r)^{1+\epsilon}\Delta u}^{\Om_{\frac\delta2}}
\\&\underset{\T{Theorem~\ref{t3.2}}}\lesssim\delta^{1-\epsilon}.
\end{split}
\end{equation}
By \eqref{3.0}, we have for a suitable $c$
$$
\tilde u_{\frac\delta2}\le u_{\frac\delta2}\le u+c\delta^{\min(\frac{\alpha}{m},\frac{\gamma}{m})}\quad\T{ in a neighborhood of $b\Om_\delta$}.
$$
We are going to  apply Theorem \ref{tse} for $\Om_\delta$ with $U:=u+c\delta^{\min(\frac{\alpha}{m},\frac{\gamma}{m})}$, $W: =\tilde u_{\frac\delta2}$ and $s:=1$; thus we get

\begin{equation}
\Label{3.7}
\begin{split}
\underset{\Om_\delta}\sup\left(\tilde u_{\frac\delta2}-(u+c\delta^{\min(\frac{\alpha}{m},\frac{\gamma}{m})})\right)&\underset{\T{stability estimate}}\lesssim \no{\left(\tilde u_{\frac\delta2}-(u+c\delta^{\min(\frac{\alpha}{m},\frac{\gamma}{m})})\right )_+}^\eta_{L^1(\Om_\delta)}
\\
&\lesssim \no{\tilde u_{\frac\delta2}-u}^\eta_{L^1(\Om_\delta)}\\
&\underset{\T{\eqref{hypernova}}}\lesssim \delta^{(1-\epsilon)\eta},
\end{split}
\end{equation}
for any $\eta<\frac{1}{2}\gamma_p=\frac{1}{np+1}$.
It follows
$$\underset{\Om_\delta}\sup\left(\tilde u_{\frac\delta2}-u\right)\lesssim\delta^{(1-\epsilon)\eta}+\delta^{\min(\frac{\alpha}{m},\frac{\gamma}{m})}\lesssim \delta^{\frac\gamma2}+ \delta^{\min(\frac{\alpha}{m},\frac{\gamma}{m})}
$$
and hence \eqref{*} is proved since $m\ge 2$. Here the last inequality follows by choosing $\epsilon=\frac{\gamma_p-\gamma}{\gamma_p+\gamma}>0$ and $\eta=\frac{1}{4}(\gamma_p+\gamma)< \frac{1}{2}\gamma_p$ since $\gamma<\gamma_p$.\\

From \eqref{3.0} and \eqref{a} (which is equivalent to \eqref{*}), it is easy to prove that 
$$|u(z)-u(z')|\lesssim |z-z'|^{\min(\frac{\alpha}{m},\frac{\gamma}{m})}\quad\T{for any $z,z'\in \bar{\Om}$};$$ thus the proof of Theorem~\ref{t1.1} is complete. 
\epf

\end{document}